# Bounds on the Number of Numerical Semigroups of a Given Genus

Maria Bras-Amorós *

October 31, 2018

**Abstract**

Combinatorics on multisets is used to deduce new upper and lower bounds on the number of numerical semigroups of each given genus, significantly improving existing ones. In particular, it is proved that the number $n_g$ of numerical semigroups of genus $g$ satisfies $2F_g \leqslant n_g \leqslant 1 + 3 \cdot 2^{g-3}$, where $F_g$ denotes the $g$th Fibonacci number.

## 1 Introduction

Let $\mathbb{N}_0$ denote the set of all non-negative integers. A *numerical semigroup* is a subset $\Lambda$ of $\mathbb{N}_0$ containing 0, closed under summation and with finite complement in $\mathbb{N}_0$. The elements in the complement $\mathbb{N}_0 \setminus \Lambda$ are called the *gaps* of the numerical semigroup and $|\mathbb{N}_0 \setminus \Lambda|$ is its *genus*. The largest gap is the *Frobenius number* of $\Lambda$ and it is at most two times the genus minus one. If it equals this bound then the numerical semigroup is said to be symmetric.

Some results have been proved related to the number of numerical semigroups of a given Frobenius number [1] and the number of symmetric semigroups of a given Frobenius number (and thus, the number of symmetric semigroups of a given genus) [4, 8]. In this work we address the problem of counting the number of numerical semigroups of a given genus.

We denote by $n_g$ the number of numerical semigroups of genus $g$. It is easy to check that $n_0 = 1$ and $n_1 = 1$. The values up to $n_{16}$ where computed by Nivaldo Medeiros and Shizuo Kakutani, and the values up to $n_{50}$ can be found in [2]. It is proved in [3] that any numerical semigroup can be represented by a unique Dyck path of order given by its genus and thus $n_g \leqslant C_g$ where $C_g$ denotes the Catalan number, $C_g = \frac{1}{g+1}\binom{2g}{g}$. It is conjectured in [2] that the sequence $(n_g)$ asymptotically behaves like the Fibonacci sequence. More precisely, $n_g \geqslant n_{g-1}+n_{g-2}$, for $g \geqslant 2$; $\lim_{g\to\infty} \frac{n_{g-1}+n_{g-2}}{n_g} = 1$; $\lim_{g\to\infty} \frac{n_g}{n_{g-1}} = \phi$, where $\phi$ is the golden ratio. See [5] for further results in this direction.

*This work was partly supported by the Spanish Ministry of Education through projects TSI2007-65406-C03-01 "E-AEGIS" and CONSOLIDER CSD2007-00004 "ARES", and by the Government of Catalonia under grant 2005 SGR 00446.

Let $F_i$ denote the $i$th Fibonacci number starting by $F_0 = 0$, $F_1 = 1$. We prove that

$$2F_g \leqslant n_g \leqslant 1 + 3 \cdot 2^{g-3}.$$

## 2 Some Results on Combinatorics

**Lemma 1.** *The multisets $A_g$ defined recursively by $A_2 = \{1, 3\}$,*

$$A_g = \{g+1\} \cup \left( \bigcup_{m \in A_{g-1}} \{0, 1, \ldots, m-1\} \right) \setminus \{g-2\}$$

*for $g > 2$ (see Figure 1) satisfy, if $g \geqslant 2$,*

$$A_g = \left( \{\overbrace{0, 0, \ldots, 0}^{2F_{g-2}}\} \cup \{\overbrace{1, 1, \ldots, 1}^{2F_{g-3}}\} \cup \{\overbrace{2, 2, \ldots, 2}^{2F_{g-4}}\} \cup \ldots \cup \{\overbrace{g-4, g-4}^{2F_2}\} \cup \{\overbrace{g-3, g-3}^{2F_1}\} \right) \cup \{g-1, g+1\}$$

*and*

$$|A_g| = 2F_g.$$

*Proof.* Both results can be proved by induction and are a consequence from the fact that, for $i \geqslant 2$, $F_i = 1 + \sum_{j=1}^{i-2} F_j$. This in turn can be proved by induction. Indeed, it is obvious for $i = 2$. If $i > 2$, by the induction hypothesis $F_{i-1} = 1 + \sum_{j=1}^{i-3} F_j$ and hence $F_i = F_{i-1} + F_{i-2} = 1 + \sum_{j=1}^{i-2} F_j$. □

$A_2 = \{1, 3\}$
$A_3 = \{0, 0, 2, 4\}$
$A_4 = \{0, 0, 1, 1, 3, 5\}$
$A_5 = \{0, 0, 0, 0, 1, 1, 2, 2, 4, 6\}$
$A_6 = \{0, 0, 0, 0, 0, 0, 1, 1, 1, 1, 2, 2, 3, 3, 5, 7\}$
$A_7 = \{0, 0, 0, 0, 0, 0, 0, 0, 0, 0, 1, 1, 1, 1, 1, 1, 2, 2, 2, 2, 3, 3, 4, 4, 6, 8\}$

Figure 1: First multisets $A_g$ as in Lemma 1.

**Lemma 2.** *The multisets $B_g$ defined recursively by $B_2 = \{1, 3\}$,*

$$B_g = \{0, g+1\} \cup \left( \bigcup_{m \in B_{i-1}} \{1, 2, \ldots, m\} \right) \setminus \{g, g-2\}$$

*for $g > 2$ (see Figure 2) satisfy, if $g > 2$,*

$$B_g = \{0\} \cup \left( \{\overbrace{1,1,\dots,1}^{3\cdot 2^{g-4}}\} \cup \{\overbrace{2,2,\dots,2}^{3\cdot 2^{g-5}}\} \cup \ldots \cup \{\overbrace{g-3,g-3,g-3}^{3\cdot 2^{0}}\} \right) \cup \{g-2, g-1, g+1\}$$

*and*

$$|B_g| = 1 + 3 \cdot 2^{g-3}.$$

*Proof.* Both results can be proved by induction and are a consequence from the fact that, for $i \geqslant 0$, $2^i = 1 + \sum_{j=0}^{i-1} 2^j$. This in turn can be proved by induction. □

$B_2 = \{1, 3\}$
$B_3 = \{0, 1, 2, 4\}$
$B_4 = \{0, 1, 1, 1, 2, 3, 5\}$
$B_5 = \{0, 1, 1, 1, 1, 1, 1, 2, 2, 2, 3, 4, 6\}$
$B_6 = \{0, 1, 1, 1, 1, 1, 1, 1, 1, 1, 1, 1, 1, 2, 2, 2, 2, 2, 2, 3, 3, 3, 4, 5, 7\}$
$B_7 = \{0, 1, 1, 1, 1, 1, 1, 1, 1, 1, 1, 1, 1, 1, 1, 1, 1, 1, 1, 1, 1, 1, 1, 1, 1,$
$2, 2, 2, 2, 2, 2, 2, 2, 2, 2, 2, 2, 3, 3, 3, 3, 3, 3, 4, 4, 4, 5, 6, 8\}$

Figure 2: First multisets $B_g$ as in Lemma 2.

# 3 Taking out Generators from a Semigroup

Every numerical semigroup can be generated by a finite set of elements and a minimal set of generators is unique (see for instance [4]). Given a numerical semigroup $\Lambda$ of genus $g$ and Frobenius number $f$, $\Lambda \cup \{f\}$ is a numerical semigroup and its genus is $g-1$. So, any numerical semigroup of genus $g$ can be obtained from a numerical semigroup of genus $g-1$ by removing one element larger than its Frobenius number. It is easy to check that when removing such an element from a numerical semigroup, the set obtained is a numerical semigroup if and only if the removed element belongs to the set of minimal generators. This gives a recursive procedure to obtain all numerical semigroups of genus $g$ from all numerical semigroups of genus $g-1$ by taking out, one by one, each generator that is larger than the Frobenius number for each numerical semigroup.

We can think of a tree whose root corresponds to the numerical semigroup $\mathbb{N}_0$, each numerical semigroup of genus $g$ is a node at distance $g$ from the root, and the children of a numerical semigroup are the numerical semigroups obtained when removing one by one each of its minimal generators which are larger than its Frobenius number. This construction was already considered in [6, 8, 7]. We depicted this tree in Figure 3. We wrote $< \lambda_{i_1}, \lambda_{i_2}, \dots, \lambda_{i_n} >$ to denote the numerical semigroup generated by $\lambda_{i_1}, \lambda_{i_2}, \dots, \lambda_{i_n}$. We used boldface letters for the minimal generators that are larger than the Frobenius number.

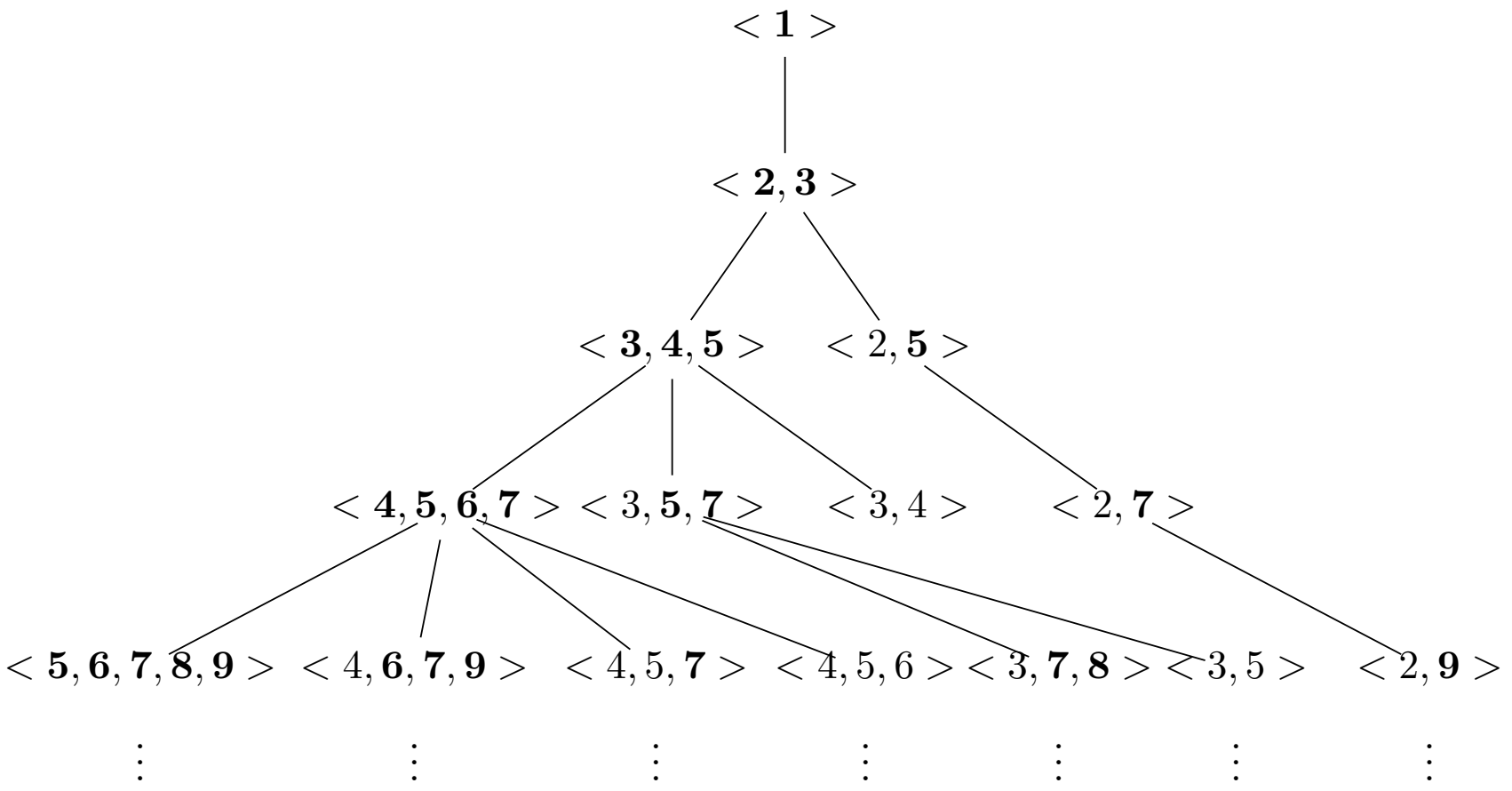

Figure 3: Recursive construction of numerical semigroups of genus $g$ from numerical semigroups of genus $g-1$.

We say that a numerical semigroup is *ordinary* if it is equal to $\{0\} \cup \{i \in \mathbb{N}_0 : i \geqslant c\}$ for some non-negative integer $c$.

**Lemma 3.** *Let $\Lambda$ be a non-ordinary numerical semigroup. Suppose that $\{\lambda_{i_1} < \lambda_{i_2} < \ldots < \lambda_{i_k}\}$ are the minimal generators of $\Lambda$ which are larger than the Frobenius number. Then the number of minimal generators of the numerical semigroup $\Lambda \setminus \{\lambda_{i_j}\}$ which are larger than its Frobenius number is*

- *at least $k-j$,*
- *at most $k-j+1$.*

*Proof.* It is obvious that the number of minimal generators which are larger than the Frobenius number is at least $k-j$ because all elements in $\Lambda \setminus \{\lambda_{i_j}\}$ which are minimal generators in $\Lambda$ are also minimal generators in $\Lambda \setminus \{\lambda_{i_j}\}$ and the new Frobenius number is $\lambda_{i_j}$.

The elements in $\Lambda \setminus \{\lambda_{i_j}\}$ which are not minimal generators in $\Lambda$ and become minimal generators in $\Lambda \setminus \{\lambda_{i_j}\}$ must be of the form $\lambda_{i_j} + \lambda_r$ for some $\lambda_r \in \Lambda$. Let $\lambda_1$ be the smallest non-zero element of $\Lambda$. If $\lambda_r > \lambda_1$ then $\lambda_{i_j} + \lambda_r - \lambda_1 > \lambda_{i_j}$ and hence $\lambda_{i_j} + \lambda_r = \lambda_1 + \lambda_s$ for some $\lambda_s \in \Lambda \setminus \{\lambda_{i_j}\}$, and $\lambda_{i_j} + \lambda_r$ is not a minimal generator of $\Lambda \setminus \{\lambda_{i_j}\}$. So, the only element that is not a minimal generator of $\Lambda$ and that may be a minimal generator of $\Lambda \setminus \{\lambda_{i_j}\}$ is $\lambda_{i_j} + \lambda_1$. $\square$

**Lemma 4.** *The ordinary semigroup $\Lambda = \{0, g+1, g+2, \ldots\}$ has minimal set of generators $\{g+1, g+2, \ldots, 2g+1\}$ and*

1. *$\Lambda \setminus \{g+1\}$ has $g+2$ minimal generators larger than its Frobenius number.*
2. *$\Lambda \setminus \{g+2\}$ has $g$ minimal generators larger than its Frobenius number.*

3. $\Lambda \setminus \{g+r\}$, *with* $r > 2$, *has* $g-r+1$ *minimal generators larger than its Frobenius number.*

*Proof.* The first item is obvious.

As proved in Lemma 3 the only element that is not a minimal generator of $\Lambda$ and that may be a minimal generator of $\Lambda \setminus \{\lambda_{i_j}\}$ is $\lambda_{i_j} + \lambda_1$. It is easy to prove that if $r = 2$ then $\lambda_{i_j} + \lambda_1$ is a minimal generator while if $r > 2$, it is not. □

**Theorem 5.** *The number* $n_g$ *of numerical semigroups of genus* $g$ *satisfies* $2F_g \leqslant n_g$ *for all* $g \geqslant 2$ *and* $2F_g \leqslant n_g \leqslant 1 + 3 \cdot 2^{g-3}$ *for all* $g \geqslant 3$.

*Proof.* Set $A_0 = B_0 = \{1\}$, $A_1 = B_1 = \{2\}$ and consider $A_g$ and $B_g$ defined as before for $g \geqslant 2$. Consider two trees $A$ and $B$ that respectively have $A_g$ and $B_g$ as the nodes at distance $g$ from its root, with the element $a \in A_g$ having $a$ children and the element $b \in B_g$ having $b$ children. It is then easy to check that, by Lemma 3 and Lemma 4, the tree in Figure 3 contains $A$ as a subtree and is contained in $B$. Thus, $|A_g| \leqslant n_g \leqslant |B_g|$. Now by Lemma 1 and Lemma 2 it follows the result. □

In Table 1 one can compare for $g$ up to 30 the actual values of $n_g$ with the bounds given in Theorem 5 and also with the bound given by the Catalan numbers proved in [3]. The values $n_g$ are from [2].

| $g$ | $2F_g$ | $n_g$ | $1 + 3 \cdot 2^{g-3}$ | $C_g$ |
|---|---|---|---|---|
| 0 | | 1 | | 1 |
| 1 | | 1 | | 1 |
| 2 | 2 | 2 | | 2 |
| 3 | 4 | 4 | 4 | 5 |
| 4 | 6 | 7 | 7 | 14 |
| 5 | 10 | 12 | 13 | 42 |
| 6 | 16 | 23 | 25 | 132 |
| 7 | 26 | 39 | 49 | 429 |
| 8 | 42 | 67 | 97 | 1430 |
| 9 | 68 | 118 | 193 | 4862 |
| 10 | 110 | 204 | 385 | 16796 |
| 11 | 178 | 343 | 769 | 58786 |
| 12 | 288 | 592 | 1537 | 208012 |
| 13 | 466 | 1001 | 3073 | 742900 |
| 14 | 754 | 1693 | 6145 | 2674440 |
| 15 | 1220 | 2857 | 12289 | 9694845 |
| 16 | 1974 | 4806 | 24577 | 35357670 |
| 17 | 3194 | 8045 | 49153 | 129644790 |
| 18 | 5168 | 13467 | 98305 | 477638700 |
| 19 | 8362 | 22464 | 196609 | 1767263190 |
| 20 | 13530 | 37396 | 393217 | 6564120420 |
| 21 | 21892 | 62194 | 786433 | 24466267020 |
| 22 | 35422 | 103246 | 1572865 | 91482563640 |
| 23 | 57314 | 170963 | 3145729 | 343059613650 |
| 24 | 92736 | 282828 | 6291457 | 1289904147324 |
| 25 | 150050 | 467224 | 12582913 | 4861946401452 |
| 26 | 242786 | 770832 | 25165825 | 18367353072152 |
| 27 | 392836 | 1270267 | 50331649 | 69533550916004 |
| 28 | 635622 | 2091030 | 100663297 | 263747951750360 |
| 29 | 1028458 | 3437839 | 201326593 | 1002242216651368 |
| 30 | 1664080 | 5646773 | 402653185 | 3814986502092304 |

Table 1: Values of $2F_g$, $n_g$, $1 + 3 \cdot 2^{g-3}$, and $C_g$ for $g$ up to 30.